\documentclass{ifacconf}

\usepackage{graphicx}      
\usepackage{natbib}        
\usepackage{subfigure} 
\usepackage{color}
\usepackage{amsmath}
\usepackage{tikz}

\newcommand{\redplus}{\tikz\draw[line width=1.5pt,red] (0,0.6ex) -- (1.2ex,0.6ex) (0.6ex,0) -- (0.6ex,1.2ex);}
\newcommand{\bluecircle}{\tikz\draw[line width=1.5pt,blue!80!white,fill=none] (0,0) circle (0.25em);}
\newcommand{\greenemptysquare}{\tikz\draw[line width=1.5pt,green!90!black,fill=none] (0,0) rectangle (0.4em,0.4em);}
\newcommand{\darkgreenemptysquare}{\tikz\draw[line width=1.5pt,black,fill=none] (0,0) rectangle (0.4em,0.4em);}
\newcommand{\redfilledcircle}{\tikz\draw[line width=1.5pt,red,fill=red] (0.2em,0.2em) circle [radius=0.15em];}
\newcommand{\purplecircle}{\tikz\draw[line width=1.5pt,red!50!blue,fill=none] (0,0) circle (0.25em);}

\begin{document}
\begin{frontmatter}

\title{An Automatic Tuning MPC with Application to Ecological Cruise Control}


\author{Mohammad Abtahi,} 
\author{Mahdis Rabbani,} 
\author{and Shima Nazari}

\address{University of California, Davis, CA 95616 USA (e-mail: sabtahi,mrabbani,snazari@ucdavis.edu)}

\begin{abstract}                
Model predictive control (MPC) is a powerful tool for planning and controlling dynamical systems due to its capacity for handling constraints and taking advantage of preview information. Nevertheless, MPC performance is highly dependent on the choice of cost function tuning parameters. In this work, we demonstrate an approach for online automatic tuning of an MPC controller with an example application to an ecological cruise control system that saves fuel by using a preview of road grade. We solve the global fuel consumption minimization problem offline using dynamic programming and find the corresponding MPC cost function by solving the inverse optimization problem. A neural network fitted to these offline results is used to generate the desired MPC cost function weight during online operation. The effectiveness of the proposed approach is verified in simulation for different road geometries.
\end{abstract}


\end{frontmatter}

\section{Introduction}
Model Predictive Control (MPC) finds the control law by solving an optimization problem over a receding horizon at each time step. MPC is especially advantageous when a preview of information is available to the controller. Therefore, it is a suitable tool for automated driving and Advanced Driver Assistance System (ADAS) applications, as a vehicle operates in a highly dynamic environment, and a preview of varying inputs from the environment can be made available to the control system through different means, such as Vehicle-to-Vehicle (V2V) communication, Global Positioning System (GPS), and sensors such as Lidar and radar. 

The performance of MPC is highly dependent on the choice of its cost function. There is no universal and systematic method for tuning the MPC cost function \cite{garriga2010model}, but rather, this tuning is usually non-intuitive and done through trial-and-error, which can be time-consuming for an inexperienced designer, \cite{4939742}. Different efforts have been reported in the literature to automate the tuning process by using model-free algorithms, \cite{tran2014model}, and optimization methods such as Particle Swarm Optimization (PSO), \cite{suzuki2007automatic}, and Genetic Algorithm (GA), \cite{sha2022automatic}. The rapid development of machine learning techniques has also encouraged researchers to use methods such as Bayesian optimization, \cite{strozecki2021automatic, sorourifar2021data}, to auto-tune the MPC cost function. 

\begin{figure}[bh]
\begin{center}
\includegraphics[width=8.4cm]{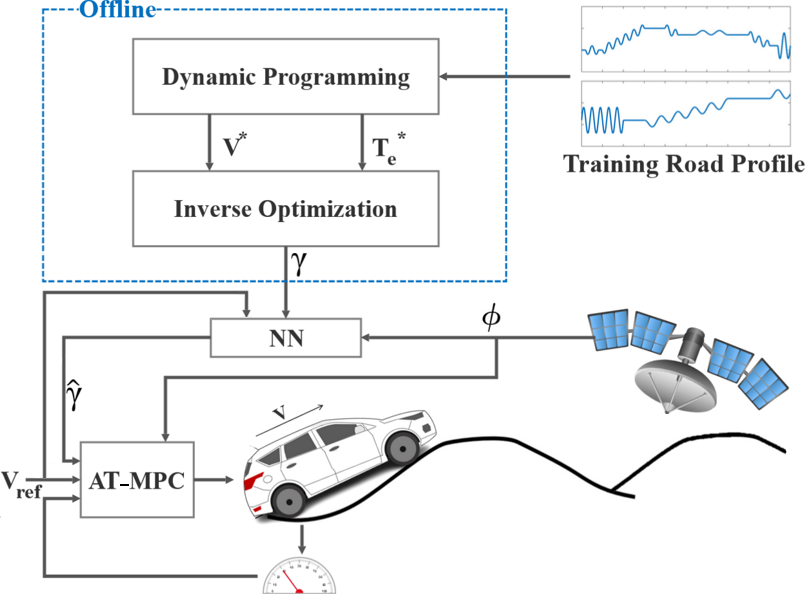}    
\caption{Summary of the Automatic Tuning MPC (AT-MPC) approach.} 
\label{fig:AT-MPC pre}
\end{center}
\end{figure}

In this work, we proposed a method for automatic tuning of an MPC controller used in an Ecological Cruise Control (ECC) system, but our approach is applicable to other systems as well. This ECC takes advantage of road elevation data to reduce a vehicle's fuel consumption. Such controller operates similarly to experienced drivers that let their vehicle accelerate over downhills to put in less gas during uphills. This strategy is especially useful for trucks and large vehicles that burn a large amount of fuel on a daily basis, \cite{hellstrom2009look}. The incorporation of road grade preview with MPC for the cruise control system is formulated in various literature, \cite{firoozi2019safe, barth2009energy}. Fuel consumption dependency on road gradient is discussed by \cite{kamal2011ecological}, where authors showed that by using the road slope information the fuel consumption can be reduced significantly. \cite{hellstrom2009look} reported up to 3.5\% fuel consumption reduction in a heavy diesel truck, and \cite{park2011predictive} showed a 4.2\% decrease in fuel consumption in a sedan car, both using optimal predictive control and road grade preview information. Note that such an optimal controller has to balance the relative importance of velocity tracking and fuel consumption at all conditions to produce good results.  Our initial study showed that finding a cost function that works well on all road elevation profiles is difficult and time-consuming, if not impossible. Therefore there is a need to systematically tune the ECC controller over different road profiles.

This work introduces a novel approach for online and automatic tuning of MPC cost function. Our approach eliminates the need for time-consuming trial-and-error in the tuning process and can adapt to varying environmental parameters and conditions during online operation, which is beyond the capability of controllers with fixed-cost functions, \cite{7782299}. The proposed method finds the global minimum fuel consumption over a road profile and solves the inverse optimization problem to find the corresponding MPC cost function weight ($\gamma$). The results from the inverse optimization and the velocity set point are used to train a Neural Network (NN) next. These steps, shown in Fig.~\ref{fig:AT-MPC pre}, are carried out offline. During online operation, the trained NN is used to tune the MPC cost function using a preview of road geometry as shown in Fig.~\ref{fig:AT-MPC pre}. We will demonstrate that our proposed method consistently produces good results on different road profiles compared to MPC controllers with a fixed cost function. The remainder of this paper is organized as follows. First, we introduce the vehicle model in Section \ref{sec:Vehicle Model}. Then the MPC formulation is presented in Section \ref{sec:ACC_LMPC}. The global fuel consumption minimization and inverse optimization problems are described in Sections \ref{sec:DP} and \ref{sec:Inverse opt}, respectively. The NN used to update the MPC cost function is presented in Section \ref{sec:NN}. Section \ref{sec:AT-MPC} describes the Auto-Tuned MPC (AT-MPC) methodology, and finally, Section \ref{sec:results} includes the results, followed by the conclusions and future research directions.

\section{Vehicle Model}
\label{sec:Vehicle Model}

In this study, we consider a Ford Escape MY 2015, which is a mid-size SUV with a 1.6 L Ecoboost engine and a six-speed automatic transmission. The full model and relevant parameters are published elsewhere, \cite{nazari2018assessing}. The longitudinal motion of the vehicle during cruising is described as
\begin{equation}\label{eq:dynamic} 
 a= \alpha_{0}T_{e}-\alpha_{1}\phi- \alpha_{2}-\alpha_{3} V- \alpha_{4} 
             {V}^2
\end{equation}
where $V$ is the vehicle velocity, and $a$ is the vehicle accelerations. The terms on the right side of the equation respectively are related to the tractive force, the gravity force, and road resistive forces (the last three terms), including rolling resistance and aerodynamic drag force. $T_e$ is the engine torque and $\phi$ is the road grade. Assuming that vehicle is in 6th gear during cruising, $\alpha_0$, $\alpha_1$, $\alpha_2$, $\alpha_3$ and $\alpha_4$ are constant coefficients computed from vehicle parameters and presented in Table~\ref{tb:parameters}.
A second-order polynomial is fitted to the vehicle fuel flow rate map as 
\begin{equation}\label{eq:fuel flow tdomain} 
    \dot m_{f}={\lambda_{0}+\lambda_{1}V+\lambda_{2}T_{e}+\lambda_{3}{T_{e}}^2+\lambda_{4}T_{e}V}+\lambda_{5}V^2
\end{equation}
where $\lambda_{i}$s are polynomial parameters and have fixed values also available in Table~\ref{tb:parameters}. The discretization domain is a controller design choice and can be made within the time domain or space domain. However, as the road information is dependent on the distance traveled by the vehicle in the studied cruise control problem, discretizing the problem in the space domain is a better choice since it eliminates the dependency of the road grade on the vehicle velocity, which is a problem state, and thus reduces the computational complexity of the problem. Equations~(\ref{eq:dynamic}) and (\ref{eq:fuel flow tdomain})  describe the change in vehicle velocity and fuel flow rate in the time domain. Equation~(\ref{eq: velocity, velocity avg}) describes these relations in the space domain, s, 
\begin{subequations}
\label{eq: velocity, velocity avg}
		\begin{align}
  \frac{dV}{ds} &=  \alpha_{0}\frac{T_{e}}{V}-\alpha_{1}\frac{\phi}{V}- \frac{\alpha_{2}}{V}-\alpha_{3}- \alpha_{4}V\\
 {\tilde {m}}_{f}= \frac{dm_f}{ds} &=  \frac{\lambda_{0}}{V}+\lambda_{1}+\lambda_{2}\frac{T_{e}}{V}+\lambda_{3}\frac{{T_{e}}^2}{V}+\lambda_{4}T_{e}+\lambda_{5} V
\end{align}
\end{subequations}
Note that the position step in our discretization is considered to be $30$ m.

\begin{table}[ht]
\begin{center}
\caption{Vehicle Model Parameters } 
\label{tb:parameters}
\begin{tabular}{c c c| c c c} 
\hline

$\lambda_{0}$& 0.5352 & $kg/h$ & $\alpha_{0}$& 0.00315 & $1/kg.m$\\
$\lambda_{1}$& -0.03021 & $kg.s/h.m$ & $\alpha_{1}$& 9.81 & $m/s$\\
$\lambda_{2}$& 0.00062 & $s^2/h.m^2$ & $\alpha_{2}$& 0.05536 & $m/s^2$\\
$\lambda_{3}$& 5.503e-5 & $s^4/kg.h .m^4$ & $\alpha_{3}$& 0.00229 & $1/s$\\
$\lambda_{4}$& 0.00079 & $s^3/h.m^3$ & $\alpha_{4}$& 2.8272e-4 & $1/m$\\
$\lambda_{5}$& 0.00131 & $kg.s^2/h.m^2$\\
\hline
\end{tabular}
\end{center}
\end{table}


\section{Cruise Controller with Linear MPC }
\label{sec:ACC_LMPC}
In this section, we introduce the Linear MPC (LMPC) problem for the cruise control system that saves fuel on hilly roads by manipulating the vehicle velocity using a preview of the road grade. The target operating conditions are moderate road slopes that do not require frequent gear shifts. Therefore, it is assumed that the vehicle will be mostly in 6th gear when cruising, and nonlinear vehicle dynamics are linearized around the cruise speed and operation in 6th gear. 
This LMPC problem can be stated as follows,
\begin{equation}\label{eq: LMPC objective}
\min\Bigg\{\gamma \bigg(\sum_{k=1}^{N} \big({\tilde {m}}_{f}^{\rm *}(k)\big)^2\bigg)+\bigg(V_{ref}-\frac{1}{N}\sum_{k=1}^{N+1} V(k)\bigg)^2\Bigg\}\\
\end{equation}\\
Subject to constraints\\
\begin{subequations}
\label{eq: LMPC const}
	\begin{align}
		\label{eq: Velocity dyn linear}
           V(k+1) &= AV(k)+B_1T_e(k) + B_2\phi(k)\\
           \label{eq:init_V}  
           V(1) &= V_i\\
		\label{eq: V_lim}
		V^{\rm min} &\leq V(k) \leq V^{\rm max}\\
       \label{eq: T_lim}
		T^{\rm min}_{e} &\leq T_{e}(k) \leq T^{\rm max}_{e}
  	\end{align}
\end{subequations}
The constraint~(\ref{eq: Velocity dyn linear}) represents the linearization of velocity dynamics in the position domain. The equality constraint~(\ref{eq:init_V}) is the initial condition, and the inequality constraints~(\ref{eq: V_lim}) and (\ref{eq: T_lim}) respectively indicate the minimum and maximum bounds on the velocity and engine torque. This controller seeks to determine the sequences of engine torque that minimize the fuel flow rate ${\tilde {m}}_{f}^{\rm *}(k)$ while keeping the average velocity over the horizon close to the reference set point $V_{ref}$. Note that in the above formulation, all velocity and engine torque values are deviations with respect to the selected linearization point, but the fuel flow is the absolute value. ${\tilde {m}}_{f}^{\rm *}(k)$ is generated by fitting a linear function to the fuel map of the engine. The variable $\gamma$ in the objective function determines the relative weight of fuel consumption against the velocity tracking and is usually selected as a fixed number in traditional MPC formulation. This trade-off between fuel consumption and velocity tracking is extremely important because poor velocity tracking discourages drivers from using this ECC system, and too tight velocity tracking undermines the fuel economy benefits of the system. 
The parameter $\gamma$ is usually selected by trial and error on various road profiles. Nevertheless, it is proven that a cost function with a fixed weight does not perform well under all operating conditions, e.g. on all road profiles, \cite{8900026}.\\ 
An example of this problem is shown in Fig.~\ref{fig: final result}, which shows fuel economy versus velocity tracking for different controllers on three different road profiles shown on the top.  Each empty circle in this figure corresponds to an LMPC controller with a different value for $\gamma$. As illustrated, some values of $\gamma$ (take $\gamma = 0.002$ as an example, shown with \purplecircle ) can produce satisfactory results on some road profiles but not on others, highlighting the necessity for online tuning of the cost function for different conditions.

\section{Global Fuel Consumption Minimization with Dynamic Programming}\label{sec:DP}

 In the first step, we formulate the fuel consumption minimization problem over the entire road profile as follows,
 \begin{equation}
\label{eq:DP cost}
\min \sum_{k=1}^{P} {\Tilde{{m}}}_{f}(k)
\end{equation} 
Subject to constraints
\begin{subequations}
\label{eq:DP const}
           \begin{align}
            \label{eq:DP vel_dynamic}
		V(k+1)&=\big(a(k)/V(k)\big)\Delta s+V(k) \\
 	    \label{eq:DP vel_avg_dynamic}
		V_{avg}(k+1) &= \frac{S(k)+\Delta s}{\frac{S(k)} 
                {V_{avg}(K)}+\frac{\Delta s}{V(K)}}
            \end{align}
            \begin{align}
		\label{eq:DP V_lim}
		V^{\rm min} &\leq V(k) \leq V^{\rm max}\\
        \label{eq:DP vavg_lim}
		V_{avg}^{\rm min} &\leq V_{avg}(k) \leq V_{avg}^{\rm max}\\
        \label{eq:Dp T_lim}
		T^{\rm min}_e &\leq T_e(k) \leq T^{\rm max}_e
            \end{align}
           \begin{align}  
        \label{eq: Dp terminal v}
                V^{\rm min} &\leq V(N) \leq  V^{\rm max}\\
        \label{eq: Dp terminal vavg}
                V_{\rm ref} &\leq V_{avg}(N) \leq  V_{avg}^{\rm max}
            \end{align}
            \begin{align}
        \label{eq: Dp initial v}
                V(1) &=V_{i}\\
        \label{eq: Dp initial vavg}
                V_{avg}(1) &=V_{i}
            \end{align}
\end{subequations}
where $P$ is the total number of steps required to travel the entire road, and ${\Tilde{{m}}}_{f}(k)$ is the fuel flow rate in the space domain. The vehicle velocity, $V(k)$, and average velocity, $V_{avg}(k)$, are problem states, shown in~(\ref{eq:DP vel_dynamic}) and (\ref{eq:DP vel_avg_dynamic}). $S(k)$ is the traveled  distance at $k$th step, and $\Delta s$ is the position discretization length. Inequality constraints (\ref{eq:DP V_lim}) and (\ref{eq:DP vavg_lim}) bound the states and (\ref{eq:Dp T_lim}) limits the input.   Inequalities (\ref{eq: Dp terminal v}) and (\ref{eq: Dp terminal vavg}) are terminal constraints. Equalities (\ref{eq: Dp initial v}) and (\ref{eq: Dp initial vavg}) indicate the initial conditions. The road grade and traveled distance are the exogenous inputs to the problem. A MATLAB-based Dynamic Programming (DP) function, \cite{sundstrom2009generic}, is used to find the global optimum solution to this problem over the different considered road profiles.
\section{Inverse Optimization Problem}\label{sec:Inverse opt} 

The solution to the discrete optimal control problem in the previous part provides the optimal trajectory of the engine torque that minimizes the fuel consumption of the vehicle while maintaining the vehicle's velocity within the predefined bounds. In the next step, we find the MPC cost function weight, $\gamma$, that generates state and input trajectories close to the global optimal solution from DP. For this, we solve the inverse optimization problem by assuming that the state and input solutions to the MPC minimization problem are known, and the only unknown is $\gamma$. Given that the MPC solution must satisfy the Karash-Kuhn Tucker (KKT) optimality condition, the value of $\gamma$ can be found. The KKT conditions are as below,
\begin{subequations}
\begin{align}
       \label{eq: KKT condition}
    \nabla f(x^*)+\sum^{m}_{i=1}p_{i}\nabla g_{i}(x^*)+\sum^{l}_{j=1} q_{j} \nabla h_{j}(x^*)=0
\end{align}
 \begin{align}
     \label{eq: KKT inequality cons1}
     q_{j}\geq 0 ,\:\:\:\:\:\;\;\;\; j= 1, ..., l
 \end{align}
  \begin{align}
     \label{eq: KKT inequality cons2}
     \sum^{l}_{j=1} q_{j} h_{j}(x^*) =0
 \end{align}
\end{subequations}
in which $x^*$ is the optimization solution, $f$ is the cost function to be minimized, $g_{i}$s are the equality constraints, $h_{j}$s are the inequality constraints, and $m$ and $l$ are the number of equality and inequality constraints, all defined as follows for the MPC problem,
\begin{subequations}
\begin{align}
       \label{eq:x^*}
    x^*=[V^*(1), ..., V^*(N+1), T^*_{e}(1), ..., T^*_{e}(N)]
\end{align}
\begin{align}
       \label{eq: KKT objective F}
   f=\gamma \bigg(\sum_{k=1}^{N} \big({\tilde {m}}_{f}^{\rm *}(k)\big)^2\bigg)+\bigg(V_{ref}-\frac{1}{N+1}\sum_{k=1}^{N+1} V(k)\bigg)^2
   \end{align}
   \begin{align}
       \label{eq: KKT equality cont}
       g_{1} &= V_{i}-V(1)\\
       g_{i+1} &= V(i+1) - AV(i)-B_1T_e(i)-B_2\phi(i) \hspace{0.1cm} i = 1,..,N\\
       h_{j} &= V^{\rm min}-V(j+1) \hspace{1.5cm} j = 1,..,N\\
       h_{j} &= V(j+1)-V^{\rm max} \hspace{1.3cm} j = N+1,..,2N\\
       h_{j} &= T^{\rm min}-T_e(j) \hspace{1.8cm} j = 2N+1,..,3N\\
       h_{j} &= T_e(j)-T_e^{\rm max} \hspace{1.75cm} j = 3N+1,..,4N
    \end{align}
\end{subequations}

We built the vector of the optimization solution, $x^*$, by concatenating vehicle velocity and engine torque for an MPC horizon of $N=60$ equal to $1.8$ km of the preview. The coefficients $p_{i}$s and $q_{j}$s are Lagrangian multipliers, which are unknown. However, $q_{j}$ are only nonzero when the inequality constraints are active (meaning that when $V(k)=V^{min}$ or $V(k)=V^{max}$, $T(k)=T^{min}$ or $T(k)=T^{max}$). Therefore, the values of $q_{j}$s are unknown only when one of the inequality constraints is active. After computing the gradient and substituting the known values, the KKT's stationary condition is simplified into,
\begin{subequations}
\begin{align}
    \label{eq:AY=b}
Q\times Y - W = 0
\end{align}
    \begin{align}
         \label{eq: Y^T}
Y^T=[\gamma, p_{1}, ..., p_{N+1}, q_{1}, ..., q_{4N}]
    \end{align}
\end{subequations}
where $Q$ and $W$ are known matrices, and $Y$ is the vector of unknowns, including the cost function weight, $\gamma$, $p_{i}$s, and unknown $q_{j}$s.  In this problem, $q_{j}$s are rarely unknown and nonzero. Thus, the above system is always an overdetermined system, meaning that the number of equations is larger than the number of unknowns. In the final step, the following optimization problem is solved to find the desired MPC cost function weight, $\gamma$,
\begin{equation}
    \label{eq: KKT norm minimization}
    \min \|\sqrt{R}QY-\sqrt{R}W\|^{2}_{2}
\end{equation}
Subject to constraints
\begin{subequations}
\begin{align}
    \label{eq: gamma const}
    \gamma &\geq 0
\end{align}
    \begin{align}
         \label{eq:qj const}
    q_{j} &\geq 0
    \end{align}
\end{subequations}
Note that the condition (\ref{eq: gamma const}) is added to ensure that the MPC cost function remains convex and thus satisfies the second-order necessary condition for optimality. R is a diagonal weight matrix for the optimization problem, and it is selected such that the terms corresponding to $V^*(1)$ and $T^*_{e}(1)$ have the weight of 1, and this weight reduces linearly to zero at the end of the MPC horizon. Using matrix $R$ simply means that we care about reproducing DP results in the near future rather than at the end of the horizon, as the MPC controller uses only the first computed control input at each step and discards the rest. It was shown in the simulations that such weighting produces a better result compared to uniform weighting.\\
The above inverse optimization problem is formulated and solved for the prediction horizon at each step for every studied road profile, and a vector of $\gamma$ is computed for each case; examples of computed $\gamma$ are shown in Fig.~\ref{fig:nnLoss} (top) and Fig.~\ref{fig:road1 predict} (bottom). Having access to the value of $\gamma$ for each part of the road, we designed an artificial neural network in the next step that can predict the desired $\gamma$ value using road features and speed set point.\\
In this work, we applied the developed methodology to four distinct road profiles. Three of these profiles are synthetic roads created by combining various sinusoidal functions to replicate the characteristics of hilly roads. Each road profile consists of several uphills, downhills, and flat sections. The maximum and minimum slopes of the profiles are 5\% and -5\%, respectively. A Test Road Profile (TRP) with a length of 270 Km is utilized for training the NN. The other two artificial Road Profiles (RP1) and (RP2) and an Actual (real-world) Road Profile (ARP) are tested to evaluate the performance of different controllers. ARP elevation data is extracted from Highway M39 NB in Michigan (US) using converted GPX files created by Google Map Pro, \cite{googlemaps}.
\section{Neural Network for MPC Cost Function Prediction}\label{sec:NN} 
Although the vector of desired $\gamma$ for the entire road profile can be determined at the beginning of each trip using DP and inverse optimization, this approach is computationally demanding, and it is more efficient to estimate $\gamma$ during online operation with a trained network.
In this paper, it is assumed that a preview of the road geometry for a longer future distance is available. A 3-km preview of the road grade ahead,  with a resolution of 30 m, and the cruise speed set point are used as inputs to the neural network to predict the desired value of $\gamma$ at each position. Note that the reason for this selection is the observed correlation between the road grade and the speed set point  with the value of $\gamma$. However, this method can be easily extended to use more features from road geometry, vehicle parameters, or environmental conditions. For the training purpose, we split the dataset, including a road profile of 270 km, cruise speed set points, and the corresponding $\gamma$ values, into testing and training sets with a portion of 20\% to 80\%, respectively, and preprocessed the data using Min-Max scaler.
The specifications of the NN are shown in Table~\ref{tb:networkProperties}. Note that the activation function is set to Rectified Linear Unit (ReLU) to avoid getting negative outputs. Furthermore, L2 regularization is deployed to avoid overfitting and improve the performance of the network for unseen inputs.
\begin{table}[hb]
\begin{center}
\caption{Neural Network Specifications}\label{tb:networkProperties}
\begin{tabular}{cc|cc}
\hline 
\textbf{Architecture} & 250-80-16 & \textbf{L2 Regularization} & $ 10^{-5} $ \\
\textbf{Optimizer} & SGD & \textbf{Loss Function} & MSE \\

\textbf{Activation} & all ReLU & \textbf{Split Portion} & 5\% \\
\textbf{Metrics} & MSE and MAE \\

\hline
\end{tabular}
\end{center}
\end{table}

\begin{figure}
\begin{center}
\includegraphics[width=8.4cm]{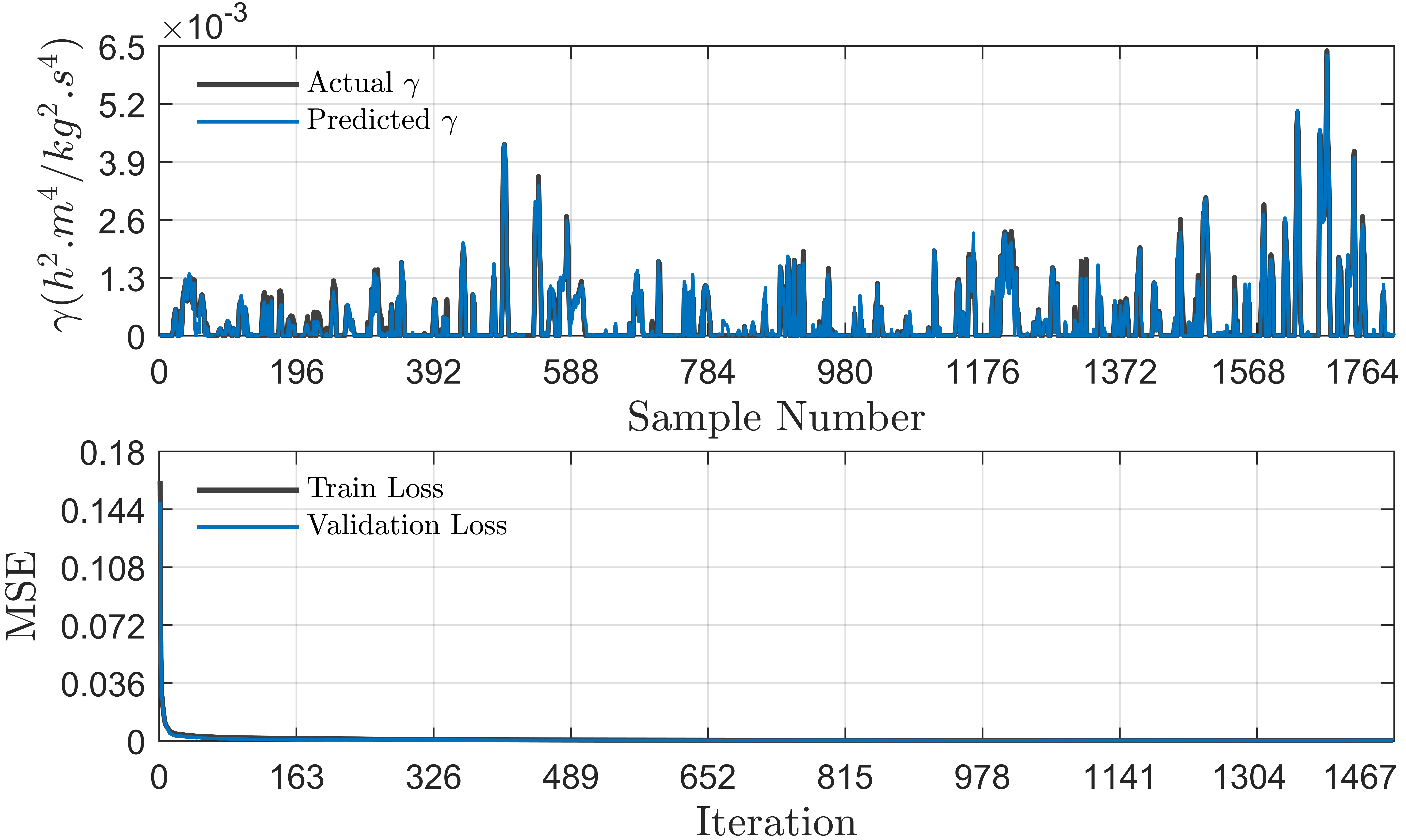}    
\caption{Training results of the neural network. (top) Actual and predicted $\gamma$ values, with the actual $\gamma$ calculated from~(\ref{eq: KKT norm minimization}). (bottom) The loss function on training and validation sets.} 
\label{fig:nnLoss}
\end{center}
\end{figure}

Fig.~\ref{fig:nnLoss} illustrates the training results of the neural network using the mentioned dataset. The top plot represents the $\gamma$ computed from~(\ref{eq: KKT norm minimization}), as the actual value of $\gamma$, and the predicted $\gamma$ based on the unseen test data. The bottom plot shows the loss function for training and validation sets. The mean squared error and mean absolute error of the test data are $3.03\times10^{-4}$ and $7.90\times10^{-3}$ on the scaled data or $1.647\times10^{-8}$ and $5.79\times10^{-5}$ on the original data, respectively.
\section{Auto-tuned MPC for ecological cruise control}\label{sec:AT-MPC} 
Using the NN, the MPC cost function can be adjusted at each step for improved and more consistent performance of the cruise controller, close to the global optimal solution from DP. As shown in Fig.~\ref{fig:AT-MPC pre}, in each iteration MPC receives the value of $\gamma$ predicted by the NN block based on a 3-km preview of the road grade ahead and the speed set point.
After updating the cost function, the proposed Auto-Tuned MPC (AT-MPC) solves the optimization problem (\ref{eq: LMPC objective}) subject to the constraints (\ref{eq: Velocity dyn linear}) to (\ref{eq: T_lim}) over a receding horizon, then it implements the first computed control input and discards the rest. This process is repeated at every step. The hard constraints on vehicle velocity are converted to soft constraints for practical reasons, and the IPOPT optimizer package is used to solve the optimization problem. The median run time for each step was $0.08$ s in the MATLAB environment. 
The controllers are tested on a nonlinear vehicle model developed in MATLAB Simulink. The model includes the vehicle longitudinal dynamics, the driveline model, including a six-speed automatic transmission model, the clutch and torque converter model, and the powertrain model, including engine torque dynamics and crankshaft dynamics. The details of the model are presented elsewhere, \cite{nazari2018assessing}.
\section{Results}\label{sec:results}
 \begin{figure}
\begin{center}
\includegraphics[width=7.9cm]{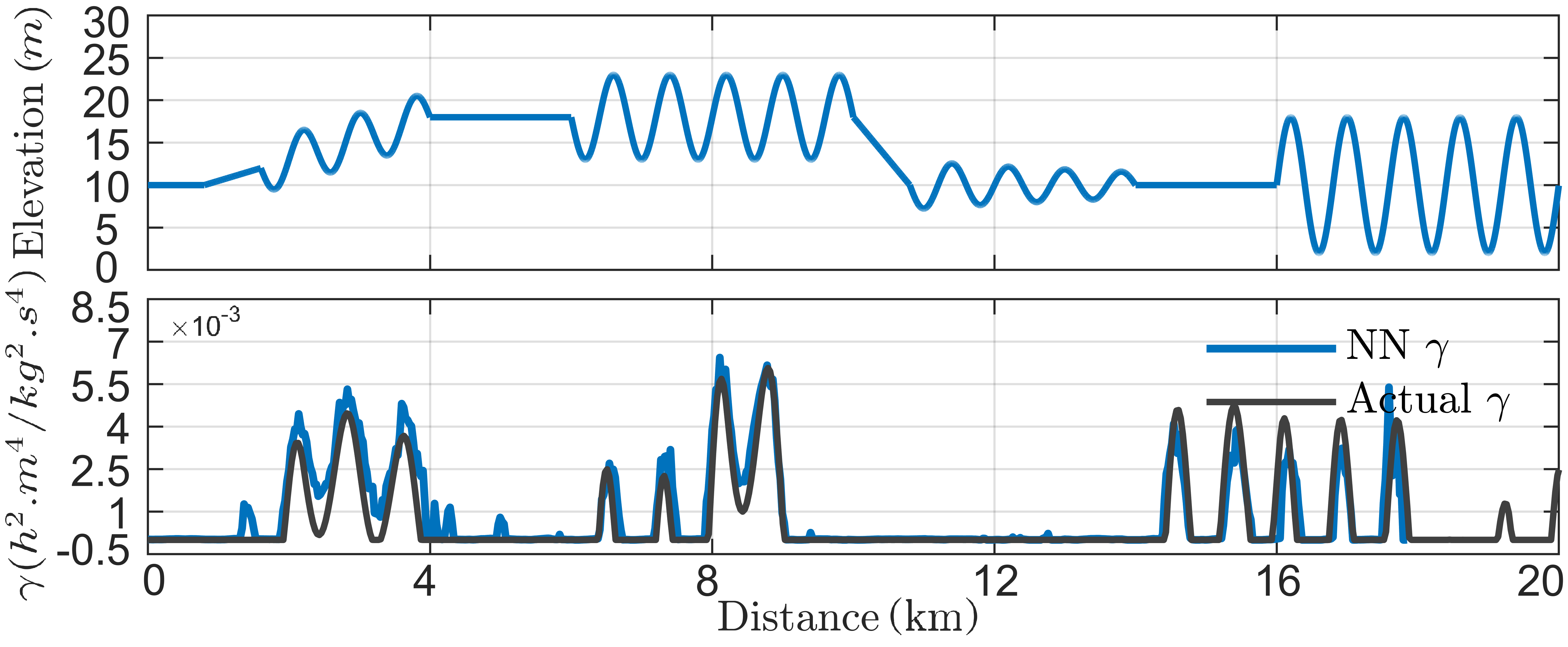}    
\caption{(top) RP1 road elevation preview, and (bottom) MPC cost function weight, $\gamma$, from inverse optimization (Actual) and predicted $\gamma$ using NN.}
\label{fig:road1 predict}
\end{center}
\end{figure}
The top plot on Fig.~\ref{fig:road1 predict} shows the road elevation for RP1 and the bottom plot compares the predicted $\gamma$ from NN with $\gamma$ from the inverse optimization results, considered as the actual value. As mentioned our target application is road profiles with moderate grades, so the reference cruise speed can be followed without gearshifts. As seen, the NN can predict the value of $\gamma$ with good accuracy on this unseen road. The accuracy of NN prediction can be improved in future work by using more data and more features in the network training.
\begin{figure}
\begin{center}
\includegraphics[width=7.9cm]{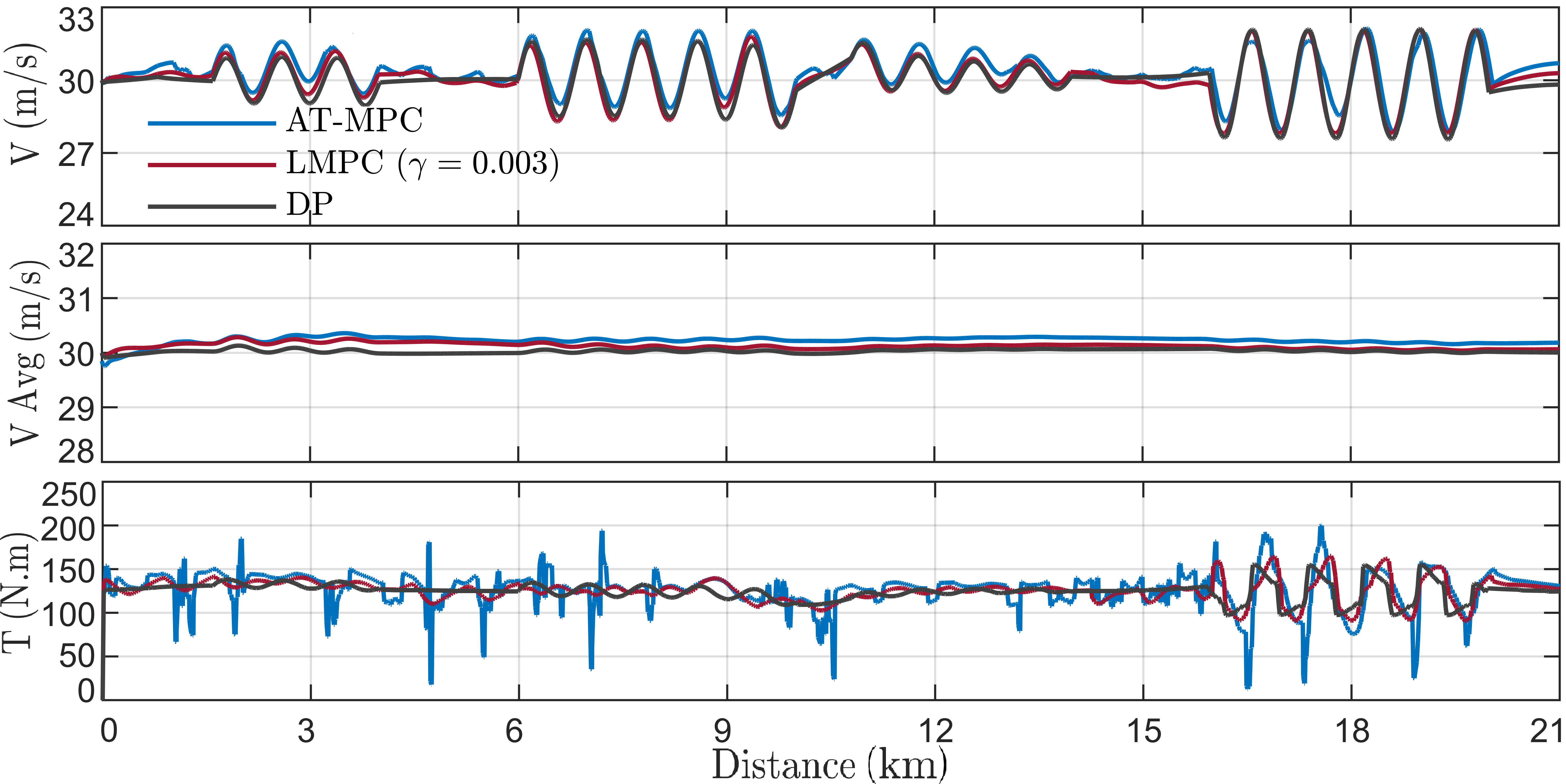}    
\caption{Comparison of the results for AT-MPC controller, LMPC with $\gamma=0.003$, and global optimal solution from DP on RP1. (top) Velocity, (middle) Velocity Average, and (bottom) Engine Torque.}
\label{fig:LPM,AT-MPC,DP}
\end{center}
\end{figure}
\begin{figure*}
\begin{center}
\includegraphics[width=17.4cm]{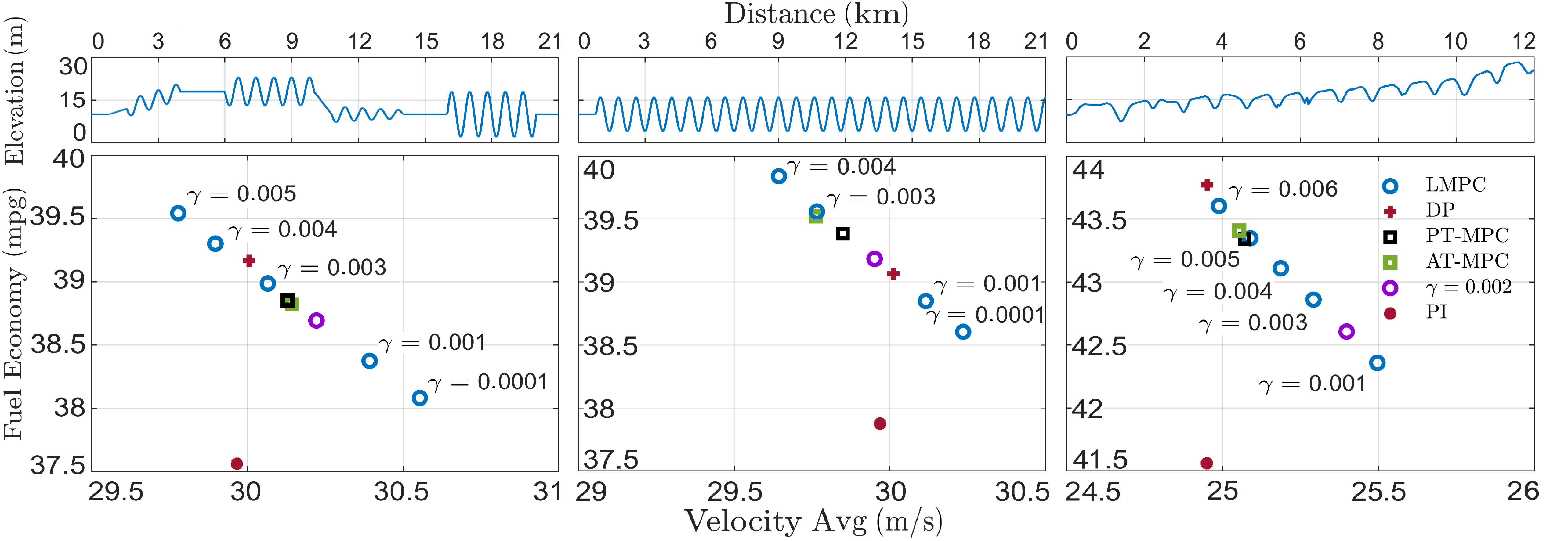}    
\caption{Fuel economy versus average velocity for DP, AT-MPC, PT-MPC, LMPC, and PI controllers for three road profiles (from left to right): RP1, RP2, and ARP.} 
\label{fig: final result}
\end{center}
\end{figure*}
Fig.~\ref{fig:LPM,AT-MPC,DP} compares the AT-MPC with the LMPC setting $\gamma= 0.003$ and the DP results for RP1 shown in Fig.~\ref{fig:road1 predict} (top). The LMPC and AT-MPC are in a closed loop with the vehicle's nonlinear Simulink model described before, and the optimal engine torque from DP is fed into the vehicle Simulink model in an open loop manner to generate the results shown in this figure and in the remainder of this paper. The $\gamma=0.003$ is selected for the LMPC as the best fixed $\gamma$ value that produces results close to DP for this road profile. The top subplot shows the vehicle velocity. The reference velocity was $30$ m/s in this simulation. As illustrated, both the LMPC and AT-MPC produce the same velocity fluctuations as in DP. The second subplot compares the average velocity. As seen, all three controllers are able to maintain an average velocity close to the reference velocity. However, the AT-MPC seems to have a higher average during the simulation compared to the other two controllers. It can indicate that AT-MPC is able to adjust its cost function in a way to allow the vehicle to accelerate more before reaching the upcoming uphill.

Finally, the last subplot compares engine torques. As seen, the torque values from AT-MPC have more fluctuations compared to DP and LMPC results. One reason for such a difference is that the cost function in the AT-MPC controller changes at each step, which gives the controller the advantage of reacting more effectively to the uphills and downhills ahead, but results in higher torque fluctuation. It should be noted that the AT-MPC torque variations are more in the low torque regions. 
Note that we have assumed none of the controllers have access to manipulate the gearshift, but the transmission control module independently controls the gear number based on an existing gearshift schedule that is published in prior work, \cite{nazari2018assessing}. In this particular case, AT-MPC performs a single-gear shift during the entire road simulation.
 \\The fuel economy versus velocity tracking for different controllers is shown in Fig.~\ref{fig: final result}, from left to right corresponding to RP1, RP2, and ARP (M39 NB Michigan). The velocity set point is $30$ m/s  for RP1 and RP2, and $25$ m/s for ARP. The elevation profile of each case is shown on top. The red cross marker (\redplus) shows the DP results generated offline and fed into the full vehicle Simulink model. The empty circle markers (\bluecircle) denote LMPC results with different cost function weights. As seen, some LMPC cases produce a good combination of fuel economy and velocity tracking, but some do not, verifying that, as explained in the LMPC section, a fixed value of $\gamma$ does not produce desirable results for all road profiles and speed set points. The dark square marker (\darkgreenemptysquare) denotes the MPC results, but with $\gamma$s directly from the inverse optimization results of~(\ref{eq: KKT norm minimization}).
 This case assumes the exact value of $\gamma$ is computed at the beginning of each trip by assuming full knowledge of the entire road geometry. Although this is a less practical approach, we included it for comparison. This case is called PreTuned MPC (PT-MPC). The light green square (\greenemptysquare) represents the AT-MPC results for each road profile. Finally, the red circle (\redfilledcircle) shows a regular cruise control system that maintains an average velocity equal to the reference velocity. This regular cruise control system is realized with a simple PI controller in the model. For each road profile, it can be seen that AT-MPC, tuned LMPCs, and DP results are mostly on the same Pareto optimal lines, which validates the effectiveness of the proposed methodology. Note that it takes a significant amount of time to find the fixed value of $\gamma$ resulting in a solution near the global optimal solution for the LMPC controller. However, the AT-MPC controller is able to tune its cost function automatically and online and produce a close result to DP on all road profiles. Comparing the results of the PI controller to DP as the best solution on the Pareto optimal line, one can observe that fuel economy is improved on average by 4\%, substantiating the effectiveness of the suggested controller in fuel consumption reduction. Furthermore, the AT-MPC results are close to the PT-MPC results, which shows that the NN can successfully produce cost function weights that are close to the actual values from the inverse optimization.

\section{Conclusion}
This paper proposed a novel method for automatic tuning of MPC cost function by solving the global minimization problem offline and inverting the optimization results to find the corresponding cost function for an MPC problem. A neural network was trained based on these offline results to generate MPC cost function weights during online operation. The method was applied to an ecological cruise control system. Simulations on a higher-fidelity vehicle model demonstrated that this approach consistently produces  good results that are close to the global optimal solution. Future research direction would be to build a larger data set by solving the inverse optimization problem for a variety of road geometries and environmental conditions and exploring more features to be included within the input features of the NN for higher prediction accuracy.
\bibliography{ifacconf}         
                                                   








\end{document}